\documentclass[12pt]{article}%
\usepackage{amsfonts}
\usepackage{sw20bams}
\usepackage{amsmath}
\usepackage{amssymb}
\usepackage{graphicx}%
\providecommand{\U}[1]{\protect\rule{.1in}{.1in}}
\begin{document}

\title{Joint Probabilities within Random Permutations}
\author{Steven Finch}
\date{May 1, 2022}
\maketitle

\begin{abstract}
A celebrated analogy between prime factorizations of integers and cycle
decompositions of permutations is explored here. \ Asymptotic formulas
characterizing semismooth numbers (possessing at most several large factors)
carry over to random permutations. \ We offer a survey of practical methods
for computing relevant probabilities of a bivariate or trivariate flavor.

\end{abstract}

\footnotetext{Copyright \copyright \ 2022 by Steven R. Finch. All rights
reserved.}Let $\Lambda_{r}$ denote the length of the $r^{\text{th}}$ longest
cycle in an $n$-permutation, chosen uniformly at random. \ If the permutation
has no $r^{\text{th}}$ cycle, then its $r^{\text{th}}$ longest cycle is
defined to have length $0$. \ The case $r=1$ has attracted widespread
attention \cite{Fi1-tcs9, Fi2-tcs9}. \ We have%
\[%
\begin{array}
[c]{ccc}%
\lim\limits_{n\rightarrow\infty}\mathbb{P}\left\{  \Lambda_{1}\leq
x\,n\right\}  =\rho\left(  \dfrac{1}{x}\right)  , &  & 0<x\leq1
\end{array}
\]
where $\rho=\rho_{1}$ is Dickman's function: \
\[%
\begin{array}
[c]{ccc}%
\xi\,\rho_{1}^{\prime}(\xi)+\rho_{1}(\xi-1)=0\text{ for }\xi>1, &  & \rho
_{1}(\xi)=1\text{ for }0\leq\xi\leq1.
\end{array}
\]
Also $\rho_{1}(\xi)=0$ for $\xi<0$. \ More generally \cite{KTP-tcs9,
Fi3-tcs9},%
\[%
\begin{array}
[c]{ccc}%
\lim\limits_{n\rightarrow\infty}\mathbb{P}\left\{  \Lambda_{2}\leq
y\,n\right\}  =\rho_{2}\left(  \dfrac{1}{y}\right)  , &  & 0<y\leq\dfrac{1}%
{2};
\end{array}
\]%
\[%
\begin{array}
[c]{ccc}%
\lim\limits_{n\rightarrow\infty}\mathbb{P}\left\{  \Lambda_{3}\leq
z\,n\right\}  =\rho_{3}\left(  \dfrac{1}{z}\right)  , &  & 0<z\leq\dfrac{1}%
{3};
\end{array}
\]%
\[%
\begin{array}
[c]{ccc}%
\lim\limits_{n\rightarrow\infty}\mathbb{P}\left\{  \Lambda_{4}\leq
w\,n\right\}  =\rho_{4}\left(  \dfrac{1}{w}\right)  , &  & 0<w\leq\dfrac{1}{4}%
\end{array}
\]
where%
\[%
\begin{array}
[c]{ccc}%
\xi\,\rho_{r}^{\prime}(\xi)+\rho_{r}(\xi-1)=\rho_{r-1}(\xi-1)\text{ for }%
\xi>1, &  & \rho_{r}(\xi)=1\text{ for }0\leq\xi\leq1
\end{array}
\]
for $r=2,3,4$. \ It is known that, as $n\rightarrow\infty$, the infinite
sequence $\frac{1}{n}\left(  \Lambda_{1},\Lambda_{2},\Lambda_{3},\Lambda
_{4},\ldots\right)  $ converges to what is called the Poisson-Dirichlet
distribution with parameter $1$. \ Our interest is in the practicalities of
computing this distribution, not for infinite sequences, but merely the finite
section $\frac{1}{n}\left(  \Lambda_{1},\Lambda_{2},\Lambda_{3},\Lambda
_{4}\right)  $. \ A special case of Billingsley's formula for the
corresponding density is \cite{Bill-tcs9, Watt-tcs9, Vrsh-tcs9, ABT1-tcs9,
Kng1-tcs9}:%
\[
f_{1234}(x,y,z,w)=\dfrac{1}{x\,y\,z\,w}\,\rho\left(  \dfrac{1-x-y-z-w}%
{w}\right)  ,
\]%
\[%
\begin{array}
[c]{ccc}%
1>x>y>z>w>0, &  & x+y+z+w<1.
\end{array}
\]
More special cases include%
\[%
\begin{array}
[c]{ccccc}%
f_{123}(x,y,z)=\dfrac{1}{x\,y\,z}\,\rho\left(  \dfrac{1-x-y-z}{z}\right)  , &
& 1>x>y>z>0, &  & x+y+z<1;
\end{array}
\]%
\[%
\begin{array}
[c]{ccccc}%
f_{12}(x,y)=\dfrac{1}{x\,y}\,\rho\left(  \dfrac{1-x-y}{y}\right)  , &  &
1>x>y>0, &  & x+y<1;
\end{array}
\]%
\[%
\begin{array}
[c]{ccc}%
f_{1}(x)=\dfrac{1}{x}\,\rho\left(  \dfrac{1-x}{x}\right)  =\dfrac{d}{dx}%
\rho_{1}\left(  \dfrac{1}{x}\right)  , &  & 1>x>0;
\end{array}
\]%
\[%
\begin{array}
[c]{ccc}%
f_{2}(y)=\dfrac{d}{dy}\rho_{2}\left(  \dfrac{1}{y}\right)  , &  & \dfrac{1}%
{2}>y>0
\end{array}
\]
and likewise for $f_{3}(z)$, $f_{4}(w)$, but no compact representations for
$f_{13}(x,z)$, $f_{14}(x,w)$, $f_{23}(y,z)$ seem to be available.

For example,%
\begin{align*}
\lim\limits_{n\rightarrow\infty}\mathbb{P}\left\{  \frac{\Lambda_{1}}{n}%
\leq\frac{1}{2}\text{ \& }\frac{\Lambda_{2}}{n}\leq\frac{1}{3}\right\}   &
=\lim\limits_{n\rightarrow\infty}\mathbb{P}\left\{  \frac{\Lambda_{1}}{n}%
\leq\frac{1}{2}\right\}  -\lim\limits_{n\rightarrow\infty}\mathbb{P}\left\{
\frac{\Lambda_{1}}{n}\leq\frac{1}{2}\text{ \& }\frac{1}{3}<\frac{\Lambda_{2}%
}{n}\leq\frac{1}{2}\right\} \\
&  =%
{\displaystyle\int\limits_{0}^{1/2}}
f_{1}(x)dx-%
{\displaystyle\int\limits_{1/3}^{1/2}}
\,%
{\displaystyle\int\limits_{1/3}^{x}}
f_{12}(x,y)dy\,dx=\rho_{1}(2)-%
{\displaystyle\int\limits_{1/3}^{1/2}}
\,%
{\displaystyle\int\limits_{1/3}^{x}}
\frac{dy\,dx}{x\,y}\\
&  =\left(  1-\ln(2)\right)  -\frac{1}{2}\ln\left(  \frac{3}{2}\right)
^{2}=0.224651842493....
\end{align*}
Call this probability $A$. \ It is associated with the blue$\,\cup
\,$magenta$\,\cup\,$green trapezoid in Figure 1, i.e., the large isosceles
triangle to the left of $y=\frac{1}{2}$ with the small orange triangle
removed. \ The probability $B$ associated with the orange$\,\cup\,$brown
triangle is clearly%
\begin{align*}
\lim\limits_{n\rightarrow\infty}\mathbb{P}\left\{  \frac{\Lambda_{2}}{n}%
>\frac{1}{3}\right\}   &  =1-\rho_{2}\left(  3\right)  =-\frac{\pi^{2}}%
{12}+\frac{\ln(3)^{2}}{2}+\operatorname*{Li}\nolimits_{2}\left(  \frac{1}%
{3}\right) \\
&  =0.147220676959....
\end{align*}
Hence%
\[
\lim\limits_{n\rightarrow\infty}\mathbb{P}\left\{  \frac{\Lambda_{1}}{n}%
>\frac{1}{2}\text{ \& }\frac{\Lambda_{2}}{n}\leq\frac{1}{3}\right\}
=1-A-B=0.628127480547...
\]
which is associated with the yellow$\,\cup\,$red$\,\cup\,$cyan trapezoid,
i.e., the large isosceles triangle to the right of $y=\frac{1}{2}$ with the
small brown triangle removed. \ Such tractable symbolics (for this specific
case) tend to obscure difficult numerics (in general) when integrating, due to
an explosive singularity of $f_{12}$ at $(x,y)=(1,0)$. \ We shall devote the
rest of this paper to simple methods for computing probabilities quickly and
accurately.
\begin{figure}
[ptb]
\begin{center}
\includegraphics[
height=5.3895in,
width=5.4077in
]%
{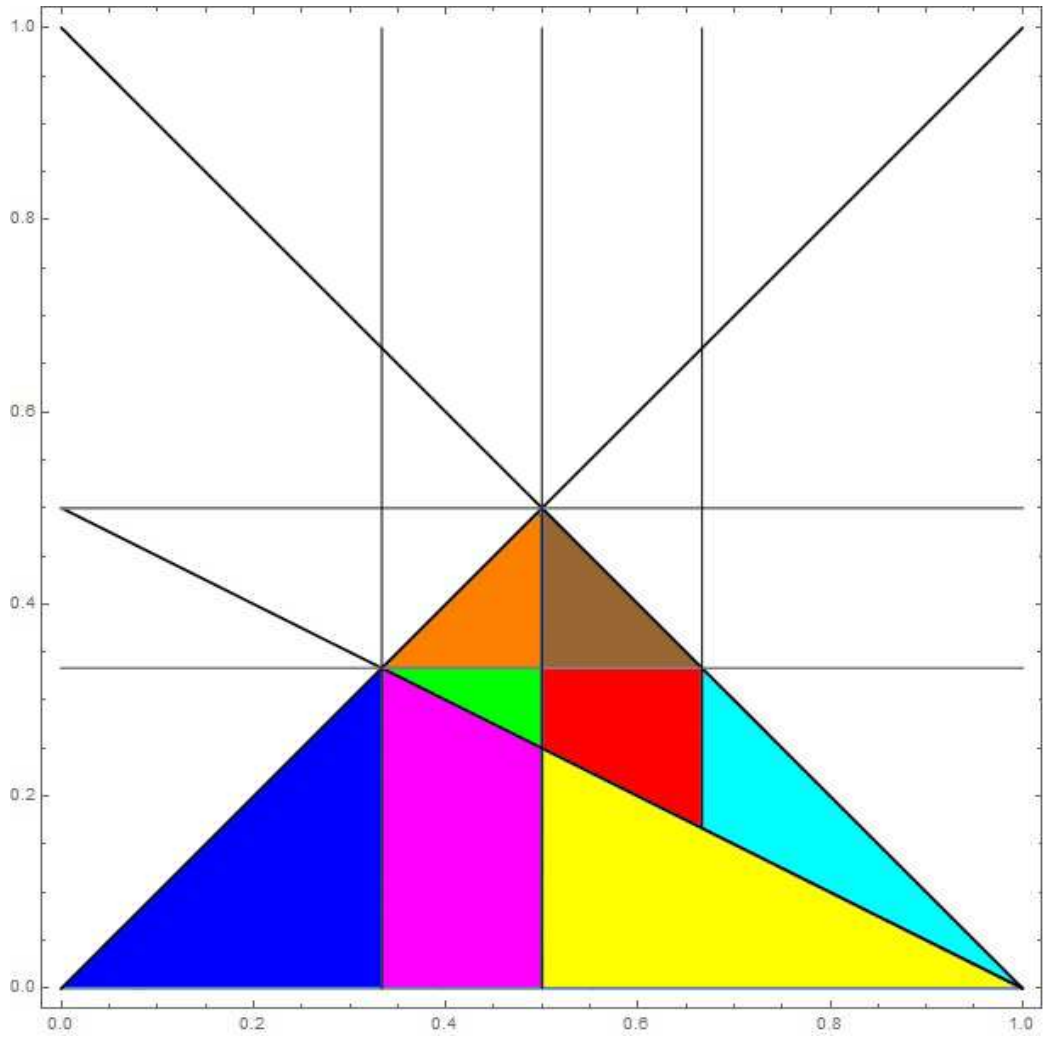}%
\caption{Domain of integration for $(\Lambda_{1},\Lambda_{2})$ example.}%
\end{center}
\end{figure}

\section{Density}

Difficulties presented by the numerical integration of $f_{12}(x,y)$ are
evident in Figure 2. \ The surface appears to touch the $xy$-plane only when
$y=0$; its prominent ridge occurs along the line $y=(1-x)/2$ because
$(1-x-y)/y=1$ corresponds to a unique point of nondifferentiability for
$\xi\mapsto\rho(\xi)$; its remaining boundary hovers over the broken line
$y=\min\{x,1-x\}$, everywhere finite except in the vicinity of $x=0$.
\begin{figure}
[ptb]
\begin{center}
\includegraphics[
height=7.3838in,
width=6.5311in
]%
{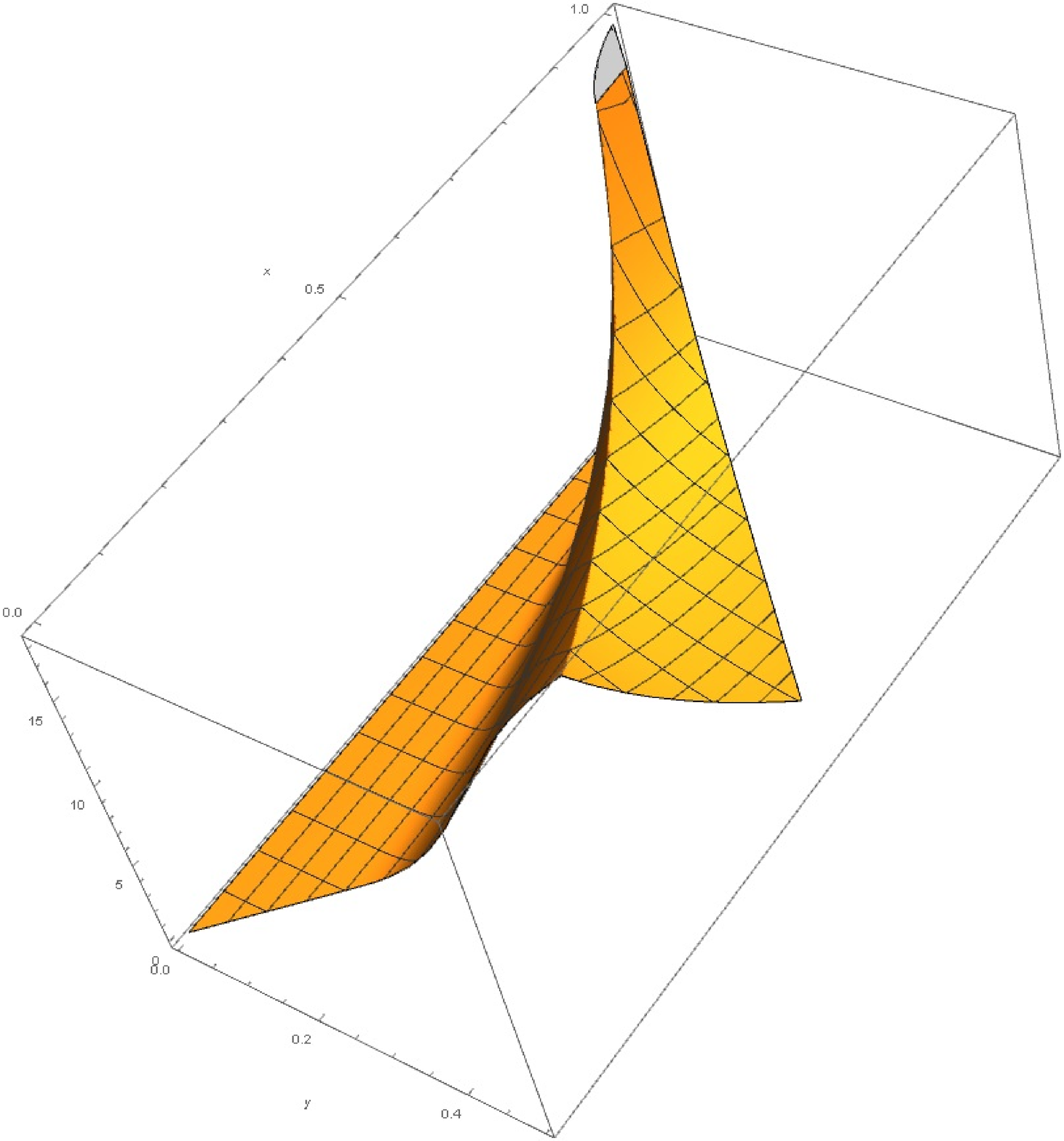}%
\caption{Probability density of $(\Lambda_{1},\Lambda_{2})$, over $0\leq
y\leq1/2$ and $y\leq x\leq1-y.$}%
\end{center}
\end{figure}

Complications are compounded for the three other densities (which are, in
themselves, approximations). \ Figure 3 contains a plot of%
\[
f_{13}(x,z)=%
{\displaystyle\int\limits_{z}^{x}}
f_{123}(x,y,z)dy.
\]
The surface appears to touch the $xz$-plane when $z=0$ and $0<x<1/2$
simultaneously, as well as everywhere along the broken line $z=\min
\{x,(1-x)/2\}$.%

\begin{figure}
[ptb]
\begin{center}
\includegraphics[
height=7.958in,
width=5.0678in
]%
{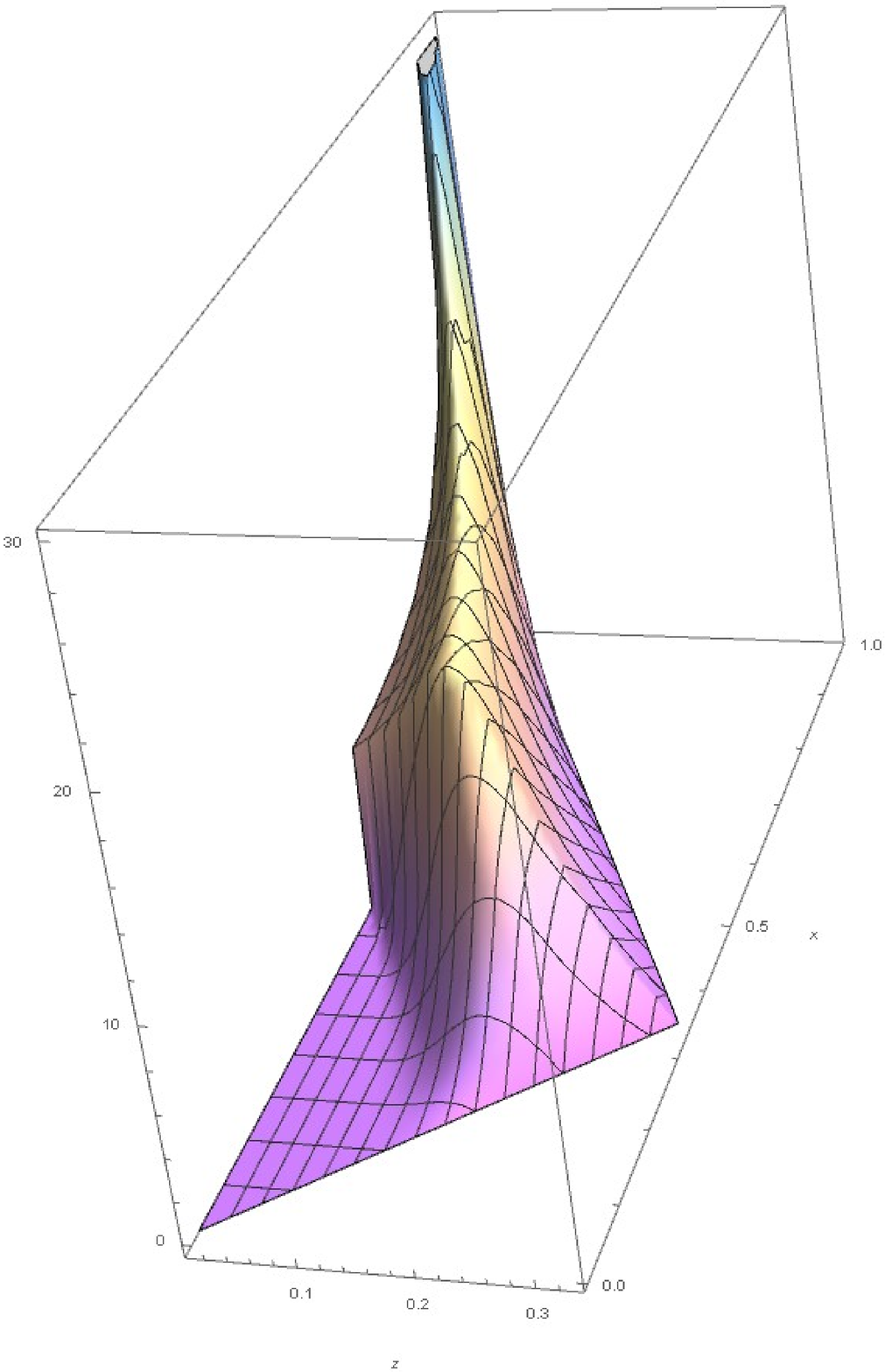}%
\caption{Probability density of $(\Lambda_{1},\Lambda_{3})$, over $0\leq
z\leq1/3$ and $z\leq x\leq1-2z.$}%
\end{center}
\end{figure}
Figure 4 contains a plot of%
\[
f_{14}(x,w)=%
{\displaystyle\int\limits_{w}^{\min\{x,1/3\}}}
\,%
{\displaystyle\int\limits_{z}^{x}}
f_{1234}(x,y,z,w)dy\,dz.
\]
The (precipitously rising) surface appears to touch the $xw$-plane only when
$w=0$ and $0<x<1/2$ simultaneously; its remaining boundary hovers over the
broken line $w=\min\{x,(1-x)/3\}$, everywhere finite except in the vicinity of
$x=0$.\ The vertical scale is more expansive here than for the other plots.\
\begin{figure}
[ptb]
\begin{center}
\includegraphics[
height=5.7242in,
width=6.8277in
]%
{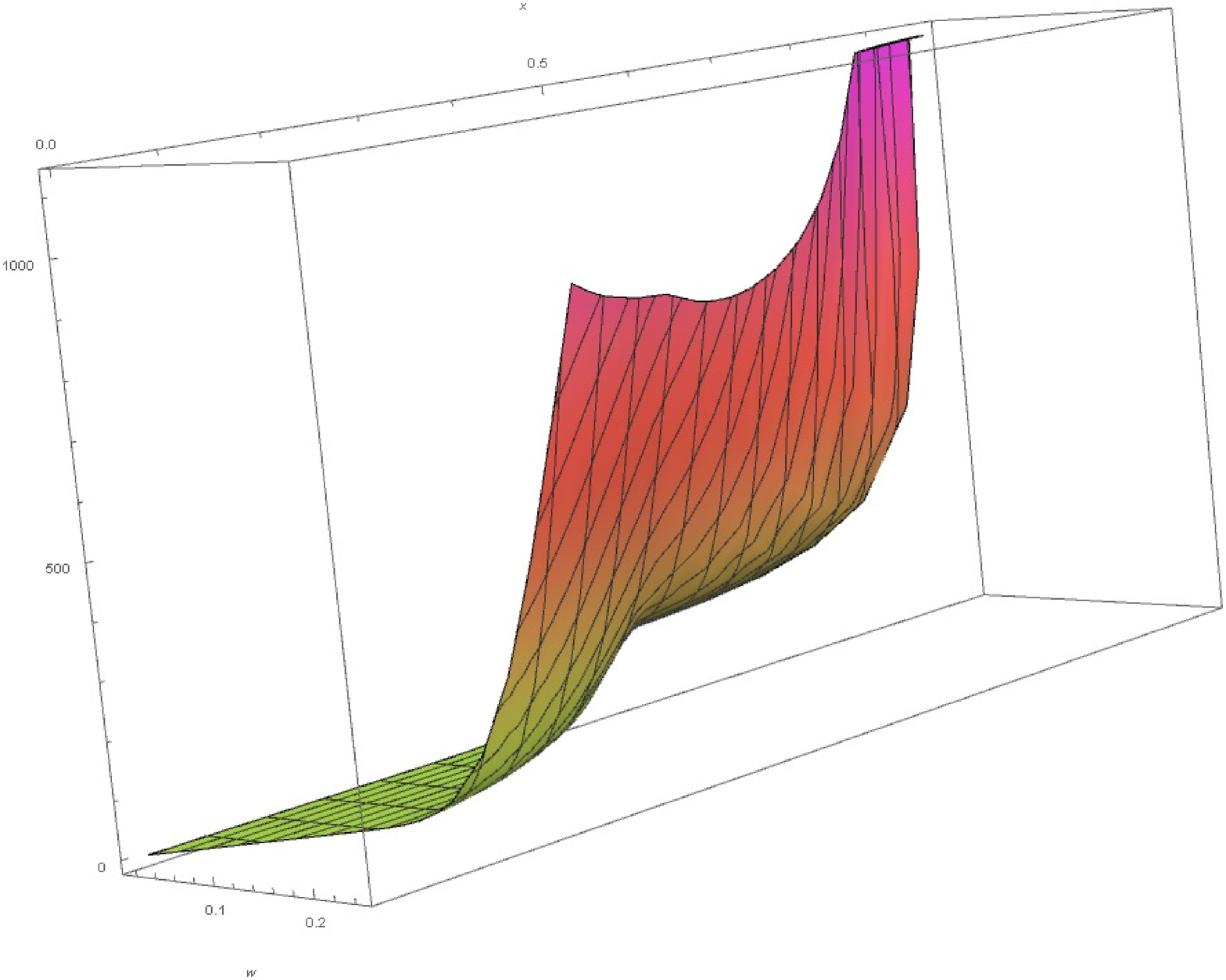}%
\caption{Probability density of $(\Lambda_{1},\Lambda_{4})$, over $0\leq
w\leq1/4$ and $w\leq x\leq1-3w.$}%
\end{center}
\end{figure}

Figure 5 contains a plot of%
\[
f_{23}(y,z)=%
{\displaystyle\int\limits_{y}^{1}}
f_{123}(x,y,z)dx.
\]
The (fairly undulating)\ surface appears to touch the $yz$-plane only when
$z=1-2y$. \ Unlike the other densities, a singularity here occurs at
$(y,z)=(0,0)$.
\begin{figure}
[ptb]
\begin{center}
\includegraphics[
height=5.8487in,
width=6.7403in
]%
{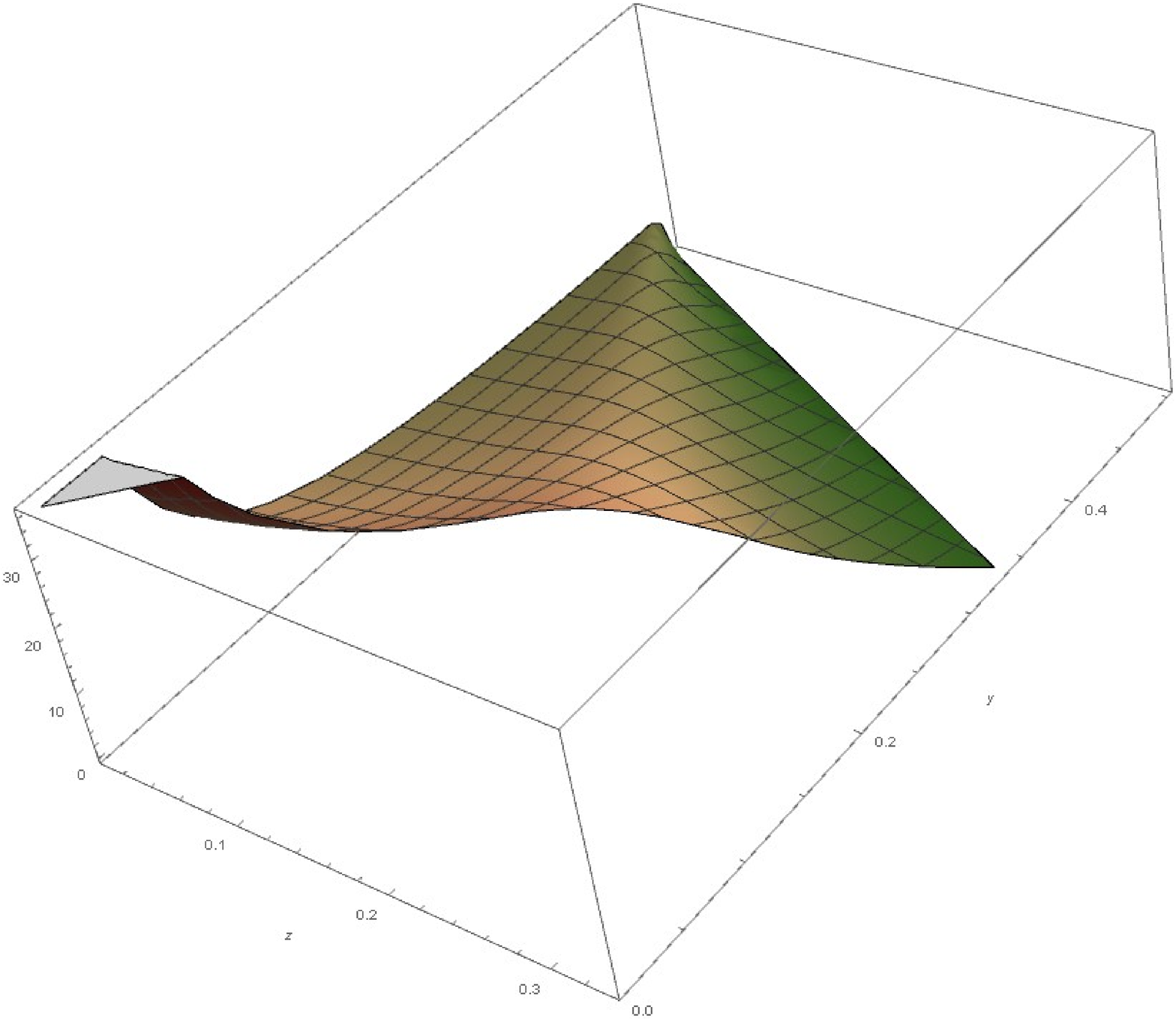}%
\caption{Probability density of $(\Lambda_{2},\Lambda_{3})$, over $0\leq
z\leq1/3$ and $z\leq y\leq(1-z)/2.$}%
\end{center}
\end{figure}

\section{Correlation}

Let%
\[%
\begin{array}
[c]{ccc}%
E(x)=%
{\displaystyle\int\limits_{x}^{\infty}}
\dfrac{e^{-t}}{t}dt=-\operatorname{Ei}(-x), &  & x>0
\end{array}
\]
be the exponential integral. \ Upon normalization, the $h^{\text{th}}$ moment
of the $r^{\text{th}}$ longest cycle length is \cite{SL-tcs9, ABT2-tcs9,
Pin-tcs9}%
\[
\lim_{n\rightarrow\infty}\frac{\mathbb{E}\left(  \Lambda_{r}^{h}\right)
}{n^{h}}=\frac{1}{h!(r-1)!}%
{\displaystyle\int\limits_{0}^{\infty}}
x^{h-1}E(x)^{r-1}\exp\left[  -E(x)-x\right]  dx
\]
(in this paper, rank $r=1,2,3$ or $4$; height $h=1$ or $2$). \ The
cross-correlation between $r^{\text{th}}$ longest and $s^{\text{th}}$ longest
cycle lengths is%
\begin{align*}
\kappa_{r,s}  &  =\frac{\mathbb{E}\left(  \Lambda_{r}\Lambda_{s}\right)
-\mathbb{E}\left(  \Lambda_{r}\right)  \mathbb{E}\left(  \Lambda_{s}\right)
}{\sqrt{\mathbb{E}\left(  \Lambda_{r}^{2}\right)  -\mathbb{E}\left(
\Lambda_{r}\right)  ^{2}}\sqrt{\mathbb{E}\left(  \Lambda_{s}^{2}\right)
-\mathbb{E}\left(  \Lambda_{s}\right)  ^{2}}}\\
&  \rightarrow\left\{
\begin{array}
[c]{lll}%
-0.75803584... &  & \text{if }r=1\text{ and }s=2,\\
-0.78421290... &  & \text{if }r=1\text{ and }s=3,\\
-0.68442819... &  & \text{if }r=1\text{ and }s=4,\\
+0.35549741... &  & \text{if }r=2\text{ and }s=3
\end{array}
\right.
\end{align*}
with cross-moments given by \cite{Grf-tcs9, Shi-tcs9}%

\[
\lim_{n\rightarrow\infty}\frac{\mathbb{E}\left(  \Lambda_{1}\Lambda
_{2}\right)  }{n^{2}}=\frac{1}{2}\,%
{\displaystyle\int\limits_{0}^{\infty}}
{\displaystyle\int\limits_{0}^{x}}
\exp\left[  -E(y)-x-y\right]  dy\,dx,
\]%
\[
\lim_{n\rightarrow\infty}\frac{\mathbb{E}\left(  \Lambda_{1}\Lambda
_{3}\right)  }{n^{2}}=\frac{1}{2}\,%
{\displaystyle\int\limits_{0}^{\infty}}
{\displaystyle\int\limits_{0}^{x}}
{\displaystyle\int\limits_{0}^{y}}
\frac{1}{y}\exp\left[  -E(z)-x-y-z\right]  dz\,dy\,dx,
\]%
\[
\lim_{n\rightarrow\infty}\frac{\mathbb{E}\left(  \Lambda_{1}\Lambda
_{4}\right)  }{n^{2}}=\frac{1}{2}\,%
{\displaystyle\int\limits_{0}^{\infty}}
{\displaystyle\int\limits_{0}^{x}}
{\displaystyle\int\limits_{0}^{y}}
{\displaystyle\int\limits_{0}^{z}}
\frac{1}{y\,z}\exp\left[  -E(w)-x-y-z-w\right]  dw\,dz\,dy\,dx,
\]%
\[
\lim_{n\rightarrow\infty}\frac{\mathbb{E}\left(  \Lambda_{2}\Lambda
_{3}\right)  }{n^{2}}=\frac{1}{2}\,%
{\displaystyle\int\limits_{0}^{\infty}}
{\displaystyle\int\limits_{0}^{x}}
{\displaystyle\int\limits_{0}^{y}}
\frac{1}{x}\exp\left[  -E(z)-x-y-z\right]  dz\,dy\,dx.
\]
The fact that $\Lambda_{1}$ is negatively correlated with other $\Lambda_{r}$,
yet $\Lambda_{2}$ is positively correlated with other $\Lambda_{s}$, is due to
longest cycles typically occupying a giant-size portion of permutations, but
second-longest cycles less so. \ 

\section{Distribution}

Bach \&\ Peralta \cite{BP-tcs9} discussed a remarkable heuristic model, based
on random bisection, that simplifies the computation of joint probabilities
involving $\Lambda_{1}$ and $\Lambda_{2}$. \ In the same paper, they
rigorously proved that asymptotic predictions emanating from the model\ are
valid. \ Subsequent researchers extended the work to $\Lambda_{1}$ and
$\Lambda_{3}$, to $\Lambda_{1}$ and $\Lambda_{4}$, and to $\Lambda_{2}$ and
$\Lambda_{3}$. \ We shall not enter into details of the model nor its absolute
confirmation, preferring instead to dwell on numerical results and certain
relative verifications.

\subsection{First and Second}

For $0<a\leq b\leq1$, Bach \&\ Peralta \cite{BP-tcs9} demonstrated that%
\[
\lim_{n\rightarrow\infty}\mathbb{P}\left\{  \frac{\Lambda_{2}}{n}\leq a\text{
\& }\frac{\Lambda_{1}}{n}\leq b\right\}  =\underset{=I_{0}(a)}{\underbrace
{\rho\left(  \frac{1}{a}\right)  }}+\underset{=I_{1}(a,b)}{\underbrace{%
{\displaystyle\int\limits_{a}^{b}}
\rho\left(  \frac{1-x}{a}\right)  \frac{dx}{x}}}.
\]
Note the slight change from earlier -- writing $\Lambda_{2}$ before
$\Lambda_{1}$ -- a convention we adopt so as to be consistent with the
literature. \ Let $J_{1}(a,b)=I_{0}(a)+I_{1}(a,b)$. \ Return now to the
example from the introduction.\ \ Evaluating%
\[
J_{1}\left(  \frac{1}{3},\frac{1}{2}\right)  =\rho(3)+%
{\displaystyle\int\limits_{1/3}^{1/2}}
\rho\left(  \frac{1-x}{1/3}\right)  \frac{dx}{x}%
\]
is less numerically problematic than evaluating%
\[
\underset{=\rho(3)}{\underbrace{%
{\displaystyle\int\limits_{0}^{1/3}}
\,%
{\displaystyle\int\limits_{0}^{x}}
f_{12}(x,y)dy\,dx}}+%
{\displaystyle\int\limits_{1/3}^{1/2}}
\,%
{\displaystyle\int\limits_{0}^{1/3}}
f_{12}(x,y)dy\,dx
\]
for two reasons:

\begin{itemize}
\item a double integral has been miraculously reduced to a single integral,

\item the argument of $\rho$ within the integral is $(1-x)/a$ rather than
$(1-x-y)/y$, which is unstable as $y\rightarrow0$.
\end{itemize}

\noindent The advantages of using the Bach \&\ Peralta formulation will become
more apparent as we move forward (incidently, their $G$ is the same as our
$J_{1}$). \ 

\begin{center}%
\begin{tabular}
[c]{|c|c||c|c|c|c|c|}\hline
$u\backslash v$ & $1$ & $1$ & $2$ & $3$ & $4$ & $5$\\\hline
$2$ & 0.30685282 & 0.69314718 &  &  &  & \\\hline
$3$ & 0.04860839 & 0.80417093 & 0.17604345 &  &  & \\\hline
$4$ & 0.00491093 & 0.61877013 & 0.09148808 & 0.01974468 &  & \\\hline
$5$ & 0.00035472 & 0.46286746 & 0.03043740 & 0.00578984 & 0.00149456 &
\\\hline
$6$ & 0.00001965 & 0.36519810 & 0.00849154 & 0.00107262 & 0.00029307 &
0.00008552\\\hline
\end{tabular}

Table 1:$\ I_{0}(1/u)$ and $I_{1}(1/u,1/v)$ for $2\leq u\leq6$, $1\leq v<u$

\bigskip%

\begin{tabular}
[c]{|c|c|c|c|c|c|c|}\hline
$u\backslash v$ & $1$ & $2$ & $3$ & $4$ & $5$ & $6$\\\hline
$2$ & 1.00000000 & 0.30685282 &  &  &  & \\\hline
$3$ & 0.85277932 & 0.22465184 & 0.04860839 &  &  & \\\hline
$4$ & 0.62368106 & 0.09639901 & 0.02465561 & 0.00491093 &  & \\\hline
$5$ & 0.46322219 & 0.03079212 & 0.00614457 & 0.00184928 & 0.00035472 &
\\\hline
$6$ & 0.36521775 & 0.00851119 & 0.00109227 & 0.00031272 & 0.00010517 &
0.00001965\\\hline
\end{tabular}

Table 2: $J_{1}(1/u,1/v)$ for $2\leq u\leq6$, $1\leq v\leq u$
\end{center}

A verification of $J_{1}(a,b)$ is as follows:%
\[
\frac{\partial J_{1}}{\partial b}=\rho\left(  \frac{1-b}{a}\right)  \frac
{1}{b}%
\]
by the Second Fundamental Theorem of Calculus, hence%
\[
\frac{\partial^{2}J_{1}}{\partial a\,\partial b}=-\rho^{\prime}\left(
\frac{1-b}{a}\right)  \frac{1-b}{a^{2}}\frac{1}{b}=\frac{\rho\left(
\dfrac{1-b}{a}-1\right)  }{\dfrac{1-b}{a}}\frac{1-b}{a^{2}b}=\frac{\rho\left(
\dfrac{1-a-b}{a}\right)  }{a\,b}=f_{12}(b,a)
\]
as anticipated by Billingsley \cite{Bill-tcs9}. \ An interpretation of
$I_{1}(a,b)$ is helpful:%
\[
I_{1}(a,b)=\lim_{n\rightarrow\infty}\mathbb{P}\left\{  \frac{\Lambda_{2}}%
{n}\leq a\text{ \& }a<\frac{\Lambda_{1}}{n}\leq b\right\}
\]
i.e., the probability that exactly one cycle has length in the interval
$(a\,n,$ $b\,n]$ and all others have length $\leq a\,n$. \ We have, for
instance,%
\[%
\begin{array}
[c]{ccc}%
\left.  \dfrac{\partial I_{1}}{\partial a}\right\vert _{b=1}=0, &  &
I_{1}(a,1)\approx0.8285
\end{array}
\]
when $a\approx0.3775\approx1/(2.649)$, the value maximizing $\mathbb{P}%
\left\{  \Lambda_{2}\leq a\,n<\Lambda_{1}\right\}  $ as $n\rightarrow\infty$. \ 

\subsection{First and Third}

For $0<a\leq1/2$ and $a\leq b\leq1$, Lambert \cite{Lmb-tcs9} demonstrated that%
\[
J_{2}(a,b)=\lim_{n\rightarrow\infty}\mathbb{P}\left\{  \frac{\Lambda_{3}}%
{n}\leq a\text{ \& }\frac{\Lambda_{1}}{n}\leq b\right\}  =J_{1}(a,b)+\underset
{=I_{2}(a,b)}{\underbrace{%
{\displaystyle\int\limits_{a}^{b}}
{\displaystyle\int\limits_{y}^{b}}
\rho\left(  \frac{1-x-y}{a}\right)  \frac{dx}{x}\frac{dy}{y}}}.
\]
(Incidently, his $G_{2}$ is the same as our $J_{2}-J_{1}=I_{2}$.)

\begin{center}%
\begin{tabular}
[c]{|c|c|c|c|c|c|}\hline
$u\backslash v$ & $1$ & $2$ & $3$ & $4$ & $5$\\\hline
$3$ & 0.14722068 & 0.08220098 &  &  & \\\hline
$4$ & 0.36143259 & 0.19556747 & 0.01998464 &  & \\\hline
$5$ & 0.46463747 & 0.20709082 & 0.02278925 & 0.00201596 & \\\hline
$6$ & 0.48588944 & 0.16644726 & 0.01263312 & 0.00136571 & 0.00013356\\\hline
\end{tabular}

Table 3: $I_{2}(1/u,1/v)$ for $3\leq u\leq6$, $1\leq v<u$

\bigskip%

\begin{tabular}
[c]{|c|c|c|c|c|c|c|}\hline
$u\backslash v$ & $1$ & $2$ & $3$ & $4$ & $5$ & $6$\\\hline
$3$ & 1.00000000 & 0.30685282 & 0.04860839 &  &  & \\\hline
$4$ & 0.98511365 & 0.29196647 & 0.04464025 & 0.00491093 &  & \\\hline
$5$ & 0.92785965 & 0.23788294 & 0.02893382 & 0.00386524 & 0.00035472 &
\\\hline
$6$ & 0.85110720 & 0.17495845 & 0.01372538 & 0.00167843 & 0.00023872 &
0.00001965\\\hline
\end{tabular}

Table 4: $J_{2}(1/u,1/v)$ for $3\leq u\leq6$, $1\leq v\leq u$
\end{center}

A verification of $J_{2}(a,b)$ is as follows:%
\begin{align*}
\frac{\partial I_{2}}{\partial b}  &  =\frac{1}{2}\frac{\partial}{\partial b}%
{\displaystyle\int\limits_{a}^{b}}
{\displaystyle\int\limits_{a}^{b}}
\rho\left(  \frac{1-x-y}{a}\right)  \frac{dx}{x}\frac{dy}{y}\\
&  =\frac{1}{2}%
{\displaystyle\int\limits_{a}^{b}}
\rho\left(  \frac{1-b-y}{a}\right)  \frac{1}{b}\frac{dy}{y}+\frac{1}{2}%
{\displaystyle\int\limits_{a}^{b}}
\rho\left(  \frac{1-x-b}{a}\right)  \frac{1}{b}\frac{dx}{x}=%
{\displaystyle\int\limits_{a}^{b}}
\rho\left(  \frac{1-x-b}{a}\right)  \frac{1}{b}\frac{dx}{x}%
\end{align*}
by symmetry; thus by Leibniz's Rule,%
\begin{align*}
\frac{\partial^{2}I_{2}}{\partial a\,\partial b}  &  =-%
{\displaystyle\int\limits_{a}^{b}}
\rho^{\prime}\left(  \frac{1-x-b}{a}\right)  \frac{1-x-b}{a^{2}}\frac{1}%
{b}\frac{dx}{x}-\rho\left(  \frac{1-a-b}{a}\right)  \frac{1}{a\,b}\\
&  =%
{\displaystyle\int\limits_{a}^{b}}
\,\frac{\rho\left(  \dfrac{1-a-x-b}{a}\right)  }{\dfrac{1-x-b}{a}}\frac
{1-x-b}{a^{2}x\,b}dx-\frac{\partial^{2}J_{1}}{\partial a\,\partial b}%
\end{align*}
hence%
\[
\frac{\partial^{2}J_{2}}{\partial a\,\partial b}=%
{\displaystyle\int\limits_{a}^{b}}
\,\frac{\rho\left(  \dfrac{1-a-x-b}{a}\right)  }{a\,x\,b}dx=%
{\displaystyle\int\limits_{a}^{b}}
f_{123}(b,x,a)dx=f_{13}(b,a),
\]
as was to be shown. \ An interpretation of $I_{2}(a,b)$ is helpful:%
\[
I_{2}(a,b)=\lim_{n\rightarrow\infty}\mathbb{P}\left\{  \frac{\Lambda_{3}}%
{n}\leq a\text{ \& }a<\frac{\Lambda_{2}}{n}\leq\frac{\Lambda_{1}}{n}\leq
b\right\}
\]
i.e., the probability that exactly two cycles have length in the interval
$(a\,n,$ $b\,n]$ and all others have length $\leq a\,n$. \ 

\subsection{First and Fourth}

For $0<a\leq1/3$ and $a\leq b\leq1$, Cavallar \cite{Cvlr-tcs9} and Zhang
\cite{Zhng-tcs9} independently demonstrated that%
\[
J_{3}(a,b)=\lim_{n\rightarrow\infty}\mathbb{P}\left\{  \frac{\Lambda_{4}}%
{n}\leq a\text{ \& }\frac{\Lambda_{1}}{n}\leq b\right\}  =J_{2}(a,b)+\underset
{=I_{3}(a,b)}{\underbrace{%
{\displaystyle\int\limits_{a}^{b}}
{\displaystyle\int\limits_{z}^{b}}
{\displaystyle\int\limits_{y}^{b}}
\rho\left(  \frac{1-x-y-z}{a}\right)  \frac{dx}{x}\frac{dy}{y}\frac{dz}{z}}}.
\]
(Incidently, Cavallar's $G_{3}$ is the same as our $J_{3}-J_{2}=I_{3}$ while
Zhang's $G_{3}$ is the same as our $J_{3}$.)

\begin{center}%
\begin{tabular}
[c]{|c|c|c|c|c|c|}\hline
$u\backslash v$ & $1$ & $2$ & $3$ & $4$ & $5$\\\hline
$4$ & 0.01488635 & 0.01488635 & 0.00396814 &  & \\\hline
$5$ & 0.07126587 & 0.06809540 & 0.01884107 & 0.00094238 & \\\hline
$6$ & 0.14082221 & 0.12382378 & 0.02870816 & 0.00222512 & 0.00009015\\\hline
\end{tabular}

Table 5: $I_{3}(1/u,1/v)$ for $4\leq u\leq6$, $1\leq v<u$

\bigskip%

\begin{tabular}
[c]{|c|c|c|c|c|c|c|}\hline
$u\backslash v$ & $1$ & $2$ & $3$ & $4$ & $5$ & $6$\\\hline
$4$ & 1.00000000 & 0.30685282 & 0.04860839 & 0.00491093 &  & \\\hline
$5$ & 0.99912552 & 0.30597834 & 0.04777489 & 0.00480762 & 0.00035472 &
\\\hline
$6$ & 0.99192941 & 0.29878222 & 0.04243355 & 0.00390355 & 0.00032887 &
0.00001965\\\hline
\end{tabular}

Table 6: $J_{3}(1/u,1/v)$ for $4\leq u\leq6$, $1\leq v\leq u$
\end{center}

We omit details of the verification of $J_{3}(a,b)$, except to mention the
start point
\[
\frac{\partial I_{3}}{\partial b}=\frac{1}{6}\frac{\partial}{\partial b}%
{\displaystyle\int\limits_{a}^{b}}
{\displaystyle\int\limits_{a}^{b}}
{\displaystyle\int\limits_{a}^{b}}
\rho\left(  \frac{1-x-y-z}{a}\right)  \frac{dx}{x}\frac{dy}{y}\frac{dz}{z}%
\]
and the end point $\partial^{2}J_{3}/\partial a\,\partial b=f_{14}(b,a)$. \ An
interpretation of $I_{3}(a,b)$ is helpful:%
\[
I_{3}(a,b)=\lim_{n\rightarrow\infty}\mathbb{P}\left\{  \frac{\Lambda_{4}}%
{n}\leq a\text{ \& }a<\frac{\Lambda_{3}}{n}\leq\frac{\Lambda_{1}}{n}\leq
b\right\}
\]
i.e., the probability that exactly three cycles have length in the interval
$(a\,n,$ $b\,n]$ and all others have length $\leq a\,n$. \ 

\subsection{Second and Third}

For $0<a<1/3$, $a\leq b<1/2$ and $b\leq c\leq1$, Ekkelkamp \cite{Ekk1-tcs9,
Ekk2-tcs9} demonstrated that%
\[
\lim_{n\rightarrow\infty}\mathbb{P}\left\{  \frac{\Lambda_{3}}{n}\leq a\text{,
}a<\frac{\Lambda_{2}}{n}\leq b\text{ \& }\frac{\Lambda_{1}}{n}\leq c\right\}
=%
{\displaystyle\int\limits_{a}^{b}}
{\displaystyle\int\limits_{y}^{c}}
\rho\left(  \frac{1-x-y}{a}\right)  \frac{dx}{x}\frac{dy}{y}%
\]
under the additional condition $a+b+c\leq1$. \ If we were to suppose that this
condition is unnecessary and set $c=1$, then by definition of $\rho_{2}$, we
would have%
\[
L_{1}(a,b)=\lim_{n\rightarrow\infty}\mathbb{P}\left\{  \frac{\Lambda_{3}}%
{n}\leq a\text{ \& }\frac{\Lambda_{2}}{n}\leq b\right\}  =\underset{=K_{0}%
(a)}{\underbrace{\rho_{2}\left(  \frac{1}{a}\right)  }}+\underset{=K_{1}%
(a,b)}{\underbrace{%
{\displaystyle\int\limits_{a}^{b}}
{\displaystyle\int\limits_{y}^{1}}
\rho_{1}\left(  \frac{1-x-y}{a}\right)  \frac{dx}{x}\frac{dy}{y}}}%
\]
where $K_{1}$ is similar (but not identical) to $I_{2}$: \
\[
K_{1}(a,b)=\lim_{n\rightarrow\infty}\mathbb{P}\left\{  \frac{\Lambda_{3}}%
{n}\leq a\text{ \& }a<\frac{\Lambda_{2}}{n}\leq b\right\}  .
\]

On the one hand, our supposition is evidently false. \ In the following, we
compare provisional theoretical values (eight digits of precision) against
simulated values (just two digits): \ 

\begin{center}%
\begin{tabular}
[c]{|l|l||l|l|l|}\hline
$u\backslash v$ & $3$ & $3$ & $4$ & $5$\\\hline
$4$ & 0.62368106 & 0.27362816 $>$ 0.21 &  & \\\hline
$5$ & 0.46322219 & 0.40043992 $>$ 0.32 & 0.17285583 $>$ 0.14 & \\\hline
$6$ & 0.36521775 & 0.43489680 $>$ 0.35 & 0.24479052 $>$ 0.20 & 0.10650591 $>$
0.09\\\hline
\end{tabular}

Table 7: $K_{0}(1/u)$ and $K_{1}(1/u,1/v)$ for $4\leq u\leq6$, $3\leq v<u$

\bigskip%

\begin{tabular}
[c]{|l|l|l|l|l|l|}\hline
$u\backslash v$ & $2$ & $3$ & $4$ & $5$ & $6$\\\hline
$3$ & 1.00000000 & 0.85277932 &  &  & \\\hline
$4$ & 0.98511365 & 0.89730922 $>$ 0.84 & 0.62368106 &  & \\\hline
$5$ & 0.92785965 & 0.86366210 $>$ 0.79 & 0.63607802 $>$ 0.60 & 0.46322219 &
\\\hline
$6$ & 0.85110720 & 0.80011455 $>$ 0.72 & 0.61000827 $>$ 0.56 & 0.47172366 $>$
0.45 & 0.36521775\\\hline
\end{tabular}

Table 8: $L_{1}(1/u,1/v)$ for $3\leq u\leq6$, $2\leq v\leq u$
\end{center}

\noindent where special cases%
\[
L_{1}(a,b)=\left\{
\begin{array}
[c]{lll}%
\rho_{2}(1/b) &  & \text{if }a=b\leq1/3,\\
\rho_{3}(1/a) &  & \text{if }a\leq1/3\text{ and }b=1/2
\end{array}
\right.
\]
are surely true.

On the other hand, a verification of $L_{1}(a,b)$ is as follows:%
\[
\frac{\partial L_{1}}{\partial b}=\frac{\partial}{\partial b}%
{\displaystyle\int\limits_{a}^{b}}
{\displaystyle\int\limits_{y}^{1}}
\rho\left(  \frac{1-x-y}{a}\right)  \frac{dx}{x}\frac{dy}{y}=%
{\displaystyle\int\limits_{b}^{1}}
\rho\left(  \frac{1-x-b}{a}\right)  \frac{1}{b}\frac{dx}{x}%
\]
hence by Leibniz's Rule,%
\begin{align*}
\frac{\partial^{2}L_{1}}{\partial a\,\partial b}  &  =-%
{\displaystyle\int\limits_{b}^{1}}
\rho^{\prime}\left(  \frac{1-x-b}{a}\right)  \frac{1-x-b}{a^{2}}\frac{1}%
{b}\frac{dx}{x}=%
{\displaystyle\int\limits_{b}^{1}}
\,\frac{\rho\left(  \dfrac{1-a-b-x}{a}\right)  }{\dfrac{1-b-x}{a}}\frac
{1-b-x}{a^{2}b\,x}dx\\
&  =%
{\displaystyle\int\limits_{b}^{1}}
\,\frac{\rho\left(  \dfrac{1-a-b-x}{a}\right)  }{a\,b\,x}dx=%
{\displaystyle\int\limits_{b}^{1}}
f_{123}(x,b,a)dx=f_{23}(b,a),
\end{align*}
as was to be shown. \ If a correction term of the form $\varphi(a)+\psi(b)$
could be incorporated into $K_{1}(a,b)$, rendering it suitably smaller, then
the above argument would still go through. \ Determining such expressions
$\varphi(a)$, $\psi(b)$ is an open problem.

For $0<\alpha<1/4$, $\alpha\leq\beta<1/3$, $\beta\leq\gamma<1/2$ and
$\gamma\leq\delta\leq1$, Ekkelkamp \cite{Ekk1-tcs9, Ekk2-tcs9} further
demonstrated that%
\begin{align*}
&  \lim_{n\rightarrow\infty}\mathbb{P}\left\{  \dfrac{\Lambda_{4}}{n}%
\leq\alpha\text{, }\alpha<\dfrac{\Lambda_{3}}{n}\leq\beta\text{, }\beta
<\dfrac{\Lambda_{2}}{n}\leq\gamma\text{ \& }\dfrac{\Lambda_{1}}{n}\leq
\delta\right\} \\
&  =%
{\displaystyle\int\limits_{\alpha}^{\beta}}
{\displaystyle\int\limits_{z}^{\gamma}}
{\displaystyle\int\limits_{y}^{\delta}}
\rho\left(  \frac{1-x-y-z}{\alpha}\right)  \frac{dx}{x}\frac{dy}{y}\frac
{dz}{x}%
\end{align*}
under the additional condition $\alpha+\beta+\gamma+\delta\leq1$. \ Such a
formula might eventually assist in calculating
\[%
\begin{array}
[c]{ccc}%
\lim\limits_{n\rightarrow\infty}\mathbb{P}\left\{  \dfrac{\Lambda_{4}}{n}%
\leq\alpha\text{ \& }\dfrac{\Lambda_{2}}{n}\leq\gamma\right\}  , &  &
\lim\limits_{n\rightarrow\infty}\mathbb{P}\left\{  \dfrac{\Lambda_{4}}{n}%
\leq\alpha\text{ \& }\dfrac{\Lambda_{3}}{n}\leq\beta\right\}  .
\end{array}
\]
We leave this task for others. \ Accuracy can be improved by including a
subordinate term -- we have studied only main terms of asymptotic expansions
-- this fact was mentioned in \cite{BS-tcs9}, citing \cite{Ekk1-tcs9}, but for
proofs one must refer to \cite{Ekk2-tcs9}. \ It is striking that so much of
this material remains unpublished (seemingly abandoned but thankfully
preserved in doctoral dissertations; see \cite{Clff-tcs9, Trmr-tcs9} for
more).\pagebreak

An odd confession is necessary at this point and it is almost surely overdue.
\ The multivariate probabilities discussed here were originally conceived not
in the context of $n$-permutations as $n\rightarrow\infty$, but instead in the
difficult realm of integers $\leq N$ (prime factorizations with cryptographic
applications) as $N\rightarrow\infty$. \ Knuth \&\ Trabb Pardo \cite{KTP-tcs9,
Grn1-tcs9, Grn2-tcs9} were the first to tenuously observe this analogy.
\ Lloyd \cite{Llyd-tcs9, Kng2-tcs9} reflected, \textquotedblleft They do not
explain the coincidence... No isomorphism of the problems is
established\textquotedblright. \ Early in his article, Tao \cite{Tao-tcs9}
wrote how a certain calculation doesn't offer understanding for
\textquotedblleft\textit{why} there is such a link\textquotedblright, but
later gave what he called a\ \textquotedblleft satisfying conceptual (as
opposed to computational) explanation\textquotedblright. \ After decades of
waiting, the fog has apparently lifted.

\section{Addendum:\ Mappings}

A counterpart of Billingsley's $f_{1234}$:%
\[
g_{1234}(x,y,z,w)=\dfrac{1}{16\,x\,y\,z\,w}\,\sigma\left(  \dfrac
{1-x-y-z-w}{w}\right)  \frac{1}{\sqrt{w}},
\]%
\[%
\begin{array}
[c]{ccc}%
1>x>y>z>w>0, &  & x+y+z+w<1;
\end{array}
\]%
\[%
\begin{array}
[c]{ccc}%
\xi\,\sigma^{\prime}(\xi)+\frac{1}{2}\sigma(\xi)+\frac{1}{2}\sigma
(\xi-1)=0\text{ for }\xi>1, &  & \sigma(\xi)=1/\sqrt{\xi}\text{ for }0<\xi
\leq1
\end{array}
\]
is applicable to the study of connected components in random mappings
\cite{Watt-tcs9, ABT1-tcs9}. \ Let $\Lambda_{1}$ and $\Lambda_{2}$ denote the
largest and second-largest such components. We use similar notation, but
different techniques (because not as much is known about $\sigma$ as about
$\rho$.) \ For example,%
\begin{align*}
\lim\limits_{n\rightarrow\infty}\mathbb{P}\left\{  \frac{\Lambda_{1}}{n}%
>\frac{1}{2}\right\}   &  =%
{\displaystyle\int\limits_{1/2}^{1}}
g_{1}(x)dx=%
{\displaystyle\int\limits_{1/2}^{1}}
\frac{1}{2x}\sigma\left(  \frac{1-x}{x}\right)  \frac{dx}{\sqrt{x}}\\
&  =\frac{1}{2}%
{\displaystyle\int\limits_{1/2}^{1}}
\frac{1}{x\sqrt{1-x}}\,dx=\ln\left(  1+\sqrt{2}\right)  .
\end{align*}
Call this probability $Q$. \ The analog here of what we called $A$ in the
introduction is%
\begin{align*}
&  1-\lim\limits_{n\rightarrow\infty}\mathbb{P}\left\{  \frac{\Lambda_{1}}%
{n}>\frac{1}{2}\right\}  -\lim\limits_{n\rightarrow\infty}\mathbb{P}\left\{
\frac{\Lambda_{1}}{n}\leq\frac{1}{2}\text{ \& }\frac{1}{3}<\frac{\Lambda_{2}%
}{n}\leq\frac{1}{2}\right\}  \\
&  =1-Q-%
{\displaystyle\int\limits_{1/3}^{1/2}}
\,%
{\displaystyle\int\limits_{1/3}^{x}}
g_{12}(x,y)dy\,dx=1-Q-%
{\displaystyle\int\limits_{1/3}^{1/2}}
\,%
{\displaystyle\int\limits_{1/3}^{x}}
\frac{1}{4\,x\,y}\sigma\left(  \frac{1-x-y}{y}\right)  \frac{dy\,dx}{\sqrt{y}%
}\\
&  =1-Q-\frac{1}{4}%
{\displaystyle\int\limits_{1/3}^{1/2}}
\,%
{\displaystyle\int\limits_{1/3}^{x}}
\frac{dy\,dx}{x\,y\sqrt{1-x-y}}=0.065484671719...
\end{align*}
and the analog of we called $1-A-B$ is%
\begin{align*}
&  \lim\limits_{n\rightarrow\infty}\mathbb{P}\left\{  \frac{\Lambda_{1}}%
{n}>\frac{1}{2}\right\}  -\lim\limits_{n\rightarrow\infty}\mathbb{P}\left\{
\frac{\Lambda_{1}}{n}>\frac{1}{2}\text{ \& }\frac{1}{3}<\frac{\Lambda_{2}}%
{n}\leq\frac{1}{2}\right\}  \\
&  =Q-%
{\displaystyle\int\limits_{1/2}^{2/3}}
\,%
{\displaystyle\int\limits_{1/3}^{1-x}}
g_{12}(x,y)dy\,dx=Q-%
{\displaystyle\int\limits_{1/2}^{2/3}}
\,%
{\displaystyle\int\limits_{1/3}^{1-x}}
\frac{1}{4\,x\,y}\sigma\left(  \frac{1-x-y}{y}\right)  \frac{dy\,dx}{\sqrt{y}%
}\\
&  =Q-\frac{1}{4}%
{\displaystyle\int\limits_{1/2}^{2/3}}
\,%
{\displaystyle\int\limits_{1/3}^{1-x}}
\frac{dy\,dx}{x\,y\sqrt{1-x-y}}=0.780087954710....
\end{align*}
Thus the analog of $B$ (associated with the orange$\,\cup\,$brown triangle in
Figure 1) is%
\[
\lim\limits_{n\rightarrow\infty}\mathbb{P}\left\{  \frac{\Lambda_{2}}{n}%
>\frac{1}{3}\right\}  =1-A-(1-A-B)=0.154427373569...
\]
and should lead in due course to a formula for $\sigma_{2}$, generalizing
$\sigma_{1}=\sigma$.%
\begin{figure}
[ptb]
\begin{center}
\includegraphics[
height=5.31in,
width=5.6498in
]%
{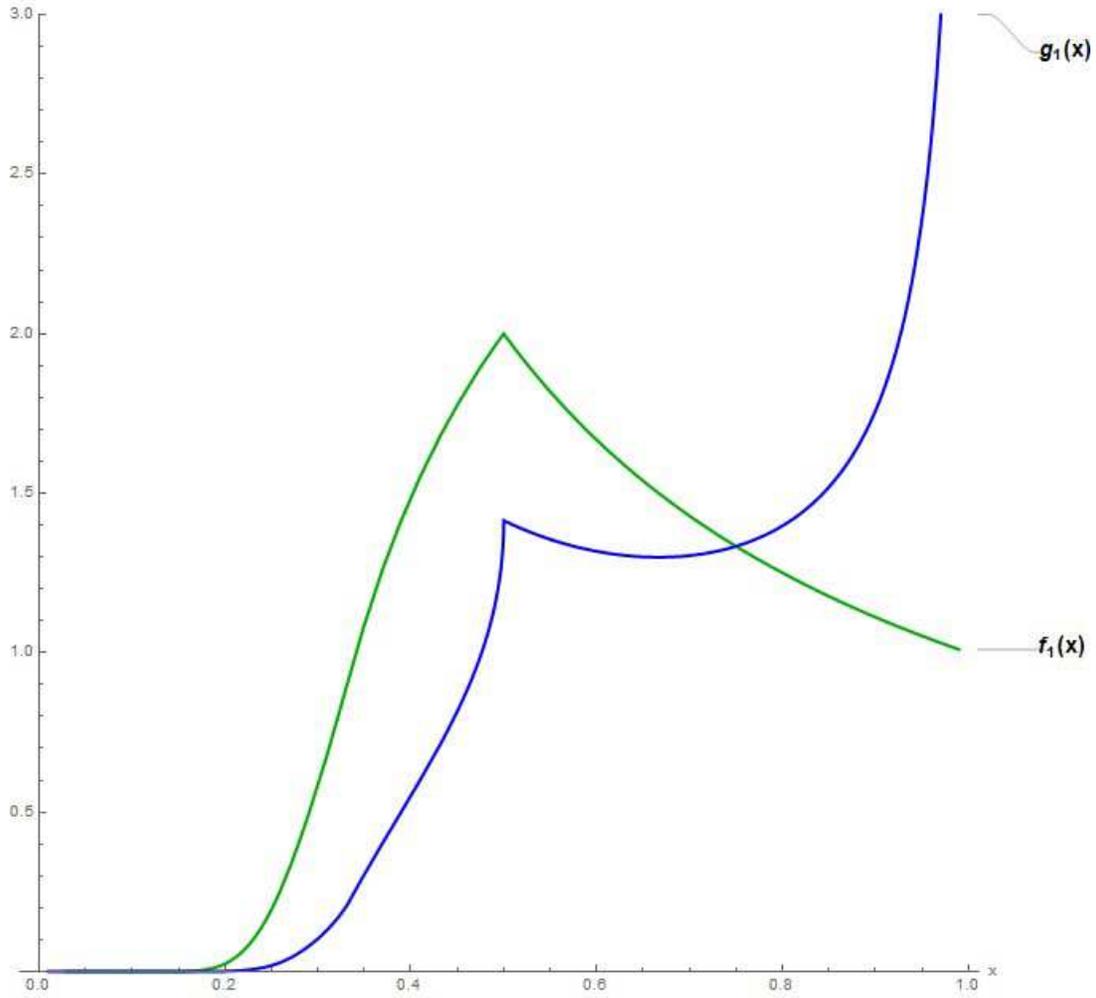}%
\caption{$f_{1}(x)=\dfrac{1}{x}\,\rho\left(  \dfrac{1-x}{x}\right)  $ and
$g_{1}(x)=\dfrac{1}{2x^{3/2}}\sigma\left(  \dfrac{1-x}{x}\right)  $
comparison;\protect\linebreak the differential expression $g_{1}(x)=\dfrac
{d}{dx}\left(  \dfrac{1}{x^{1/2}}\sigma\left(  \dfrac{1}{x}\right)  \right)  $
is akin to $f_{1}(x)=\dfrac{d}{dx}\,\rho\left(  \dfrac{1}{x}\right)  $. }%
\end{center}
\end{figure}
\begin{figure}
[ptb]
\begin{center}
\includegraphics[
height=6.5847in,
width=6.5639in
]%
{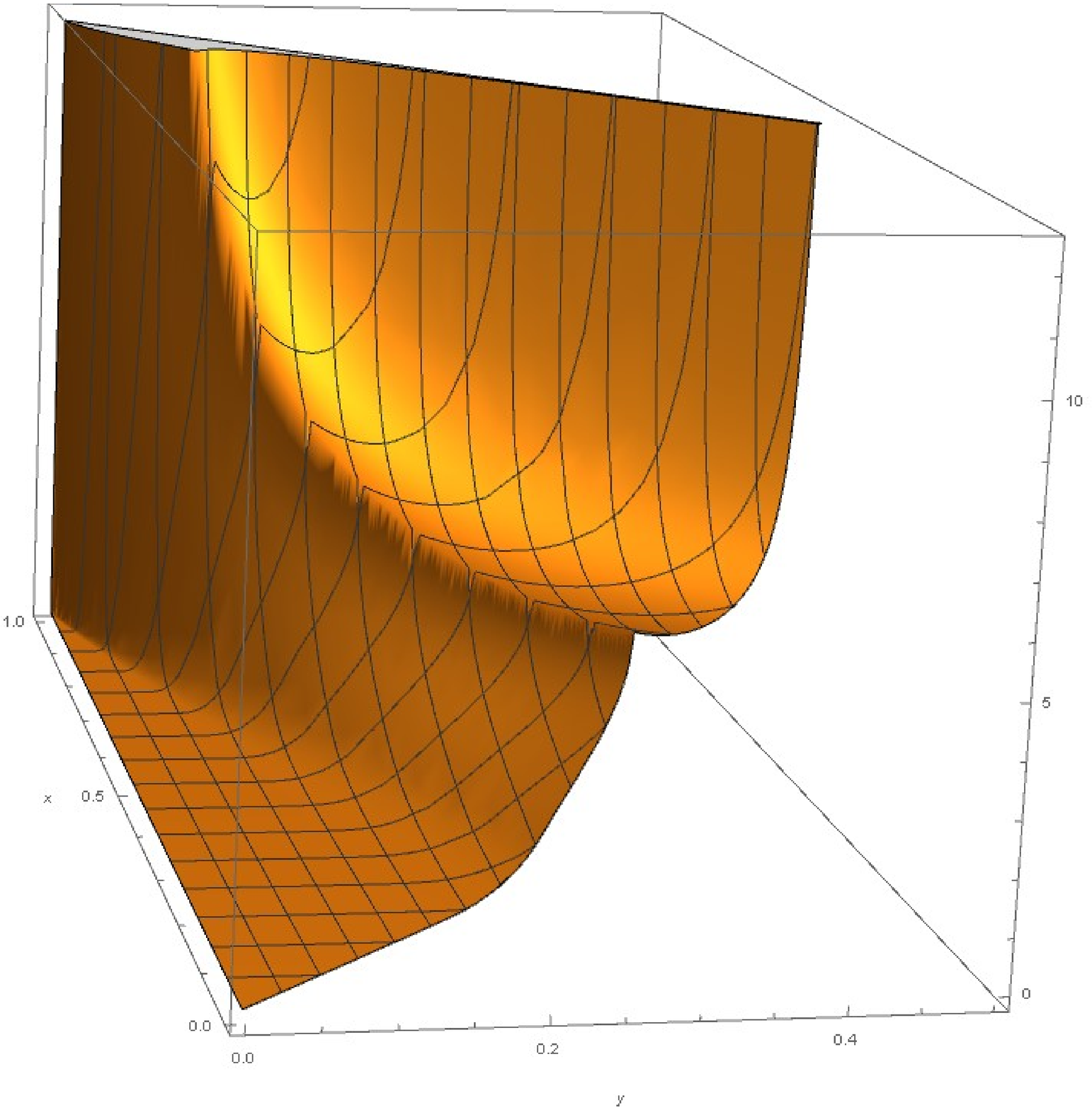}%
\caption{$g_{12}(x,y)=\dfrac{1}{4\,x\,y^{3/2}}\sigma\left(  \dfrac{1-x-y}%
{y}\right)  $ over $0\leq y\leq1/2$ and $y\leq x\leq1-y$;\medskip
\ \protect\linebreak this contrasts sharply from plot of $f_{12}(x,y)$ in
Figure 2 along diagonal segment $x=y$.}%
\end{center}
\end{figure}

\section{Addendum:\ Short Cycles}

Given a random $n$-permutation, let $S_{r}$ denote the length of the
$r^{\text{th}}$ shortest cycle ($0$ if the permutation has no $r^{\text{th}}$
cycle) and $C_{\ell}$ denote the number of cycles of length $\ell$. \ Since,
as $n\rightarrow\infty$, the distribution of $C_{\ell}$ approaches
Poisson($1/\ell$) and $C_{1}$, $C_{2}$, $C_{3}$, \ldots\ become asymptotically
independent \cite{AT-tcs9}, we can calculate corresponding probabilities for
$S_{r}$. \ For example,%
\[
\mathbb{P}\left\{  S_{1}=1\right\}  =\mathbb{P}\left\{  C_{1}\geq1\right\}
=1-\mathbb{P}\left\{  C_{1}=0\right\}  =1-e^{-1},
\]%
\begin{align*}
\mathbb{P}\left\{  S_{1}=2\right\}   &  =\mathbb{P}\left\{  C_{1}=0\text{
\&\ }C_{2}\geq1\right\}  =\mathbb{P}\left\{  C_{1}=0\right\}  -\mathbb{P}%
\left\{  C_{1}=0\text{ \&\ }C_{2}=0\right\} \\
&  =\mathbb{P}\left\{  C_{1}=0\right\}  \left(  1-\mathbb{P}\left\{
C_{2}=0\right\}  \right)  =e^{-1}\left(  1-e^{-1/2}\right)  =e^{-1}-e^{-3/2}%
\end{align*}
and, more generally,%
\[%
\begin{array}
[c]{ccc}%
\mathbb{P}\left\{  S_{1}=i\right\}  =e^{-H_{i-1}}-e^{-H_{i}}, &  & H_{m}=%
{\displaystyle\sum\limits_{k=1}^{m}}
\dfrac{1}{k}.
\end{array}
\]
It is understood that these are limiting quantities as $n\rightarrow\infty$.
\ As another example,%
\[
\mathbb{P}\left\{  S_{2}=1\right\}  =\mathbb{P}\left\{  C_{1}\geq2\right\}
=1-\mathbb{P}\left\{  C_{1}\leq1\right\}  =1-2e^{-1},
\]%
\begin{align*}
\mathbb{P}\left\{  S_{2}=2\right\}   &  =\mathbb{P}\left\{  C_{1}=1\text{
\&\ }C_{2}\geq1\right\}  +\mathbb{P}\left\{  C_{1}=0\text{ \&\ }C_{2}%
\geq2\right\} \\
&  =\mathbb{P}\left\{  C_{1}=1\right\}  -\mathbb{P}\left\{  C_{1}=1\text{
\&\ }C_{2}=0\right\}  +\mathbb{P}\left\{  C_{1}=0\right\}  -\mathbb{P}\left\{
C_{1}=0\text{ \&\ }C_{2}\leq1\right\} \\
&  =e^{-1}\left(  1-e^{-1/2}\right)  +e^{-1}\left(  1-\tfrac{3}{2}%
e^{-1/2}\right)  =2e^{-1}-\tfrac{5}{2}e^{-3/2}%
\end{align*}
and%
\[
\mathbb{P}\left\{  S_{2}=j\right\}  =\left(  H_{j-1}+1\right)  e^{-H_{j-1}%
}-\left(  H_{j}+1\right)  e^{-H_{j}}.
\]
Similar reasoning leads to
\[
\mathbb{P}\left\{  S_{1}=i\text{ \&\ }S_{2}=j\right\}  =\left\{
\begin{array}
[c]{lll}%
e^{-H_{i-1}}-\left(  1+\dfrac{1}{i}\right)  e^{-H_{i}} & \bigskip & \text{if
}i=j,\\
\dfrac{1}{i}\left(  e^{-H_{j-1}}-e^{-H_{j}}\right)  & \bigskip & \text{if
}i<j,\\
0 &  & \text{otherwise}%
\end{array}
\right.
\]
enabling a conjecture: $\mathbb{E}(S_{1}S_{2})=O(\ln(n)^{3})$. \ A\ proof
still remains out of reach.

\section{Acknowledgements}

I am grateful to Michael Rogers, Josef Meixner, Nicholas Pippenger, Eran
Tromer, John Kingman, Andrew Barbour, Ross Maller and Joseph Blitzstein for
helpful discussions. \ The creators of Mathematica, as well as administrators
of the MIT Engaging Cluster, earn my gratitude every day. \ Interest in this
subject has, for me, spanned many years \cite{Fi4-tcs9, Fi5-tcs9}. \ A\ sequel
to this paper will be released soon \cite{Fi6-tcs9}.


\begin{thebibliography}{99}                                                                                               %


\bibitem {Fi1-tcs9}S. R. Finch, Permute, Graph, Map, Derange, arXiv:2111.05720.

\bibitem {Fi2-tcs9}S. R. Finch, Rounds, Color, Parity, Squares, arXiv:2111.14487.

\bibitem {KTP-tcs9}D. E. Knuth and L. Trabb Pardo, Analysis of a simple
factorization algorithm, \textit{Theoret. Comput. Sci.} 3 (1976) 321--348;
also in \textit{Selected Papers on Analysis of Algorithms}, CSLI, 2000, pp.
303-339; MR0498355.

\bibitem {Fi3-tcs9}S. R. Finch, Second best, Third worst, Fourth in line, arXiv:2202.07621.

\bibitem {Bill-tcs9}P. Billingsley, On the distribution of large prime
divisors, \textit{Period. Math. Hungar.} 2 (1972) 283--289; MR0335462.

\bibitem {Watt-tcs9}G. A. Watterson, The stationary distribution of the
infinitely-many neutral alleles diffusion model, \textit{J. Appl. Probab.} 13
(1976) 639--651; 14 (1977) 897; MR0504014 and MR0504015.

\bibitem {Vrsh-tcs9}A. M. Vershik, Asymptotic distribution of factorizations
of natural numbers into prime divisors (in Russian), \textit{Dokl. Akad. Nauk
SSSR} v. 289 (1986) n. 2, 269--272; Engl. transl. in \textit{Soviet Math.
Dokl.} v. 34 (1987)\ 57--61; MR0856456.

\bibitem {ABT1-tcs9}R. Arratia, A. D. Barbour and S. Tavar\'{e}, Random
combinatorial structures and prime factorizations, \textit{Notices Amer. Math.
Soc.} 44 (1997) 903--910; MR1467654.

\bibitem {Kng1-tcs9}J. F. C. Kingman, Poisson processes revisited,
\textit{Probab. Math. Statist.} 26 (2006) 77--95; MR2301889.

\bibitem {SL-tcs9}L. A. Shepp and S. P. Lloyd, Ordered cycle lengths in a
random permutation, \textit{Trans. Amer. Math. Soc.} 121 (1966) 340--357; MR0195117.

\bibitem {ABT2-tcs9}R. Arratia, A. D. Barbour and S. Tavar\'{e},
\textit{Logarithmic Combinatorial Structures: a Probabilistic Approach},
Europ. Math. Society, 2003, pp. 21-24, 52, 87--89, 118; MR2032426.

\bibitem {Pin-tcs9}R. G. Pinsky, A view from the bridge spanning combinatorics
and probability, arXiv:2105.13834.

\bibitem {Grf-tcs9}R. C. Griffiths, On the distribution of allele frequencies
in a diffusion model, \textit{Theoret. Population Biol. }15 (1979) 140--158; MR0528914.

\bibitem {Shi-tcs9}T. Shi, Cycle lengths of $\theta$-biased random
permutations, B.S. thesis, Harvey Mudd College, 2014, http://scholarship.claremont.edu/hmc\_theses/65/.

\bibitem {BP-tcs9}E. Bach and R. Peralta, Asymptotic semismoothness
probabilities, \textit{Math. Comp.} 65 (1996) 1701--1715; MR1370848.

\bibitem {Lmb-tcs9}R. Lambert, \textit{Computational Aspects of Discrete
Logarithms}, Ph.D.\ thesis, Univ. of Waterloo, 1996.

\bibitem {Cvlr-tcs9}S. H. Cavallar, \textit{On the Number Field Sieve Integer
Factorisation Algorithm}, Ph.D.\ thesis, Univ. Leiden, 2002; ch. 2 also in
\textit{The Three-Large-Primes Variant of the Number Field Sieve}, CWI report
MAS-R0219, 2002, http://ir.cwi.nl/pub/4222.

\bibitem {Zhng-tcs9}C. Zhang, \textit{An Extension of the Dickman Function and
its Application}, Ph.D. thesis, Purdue Univ., 2002; Distribution of
$k$-semismooth integers, \textit{PanAmer. Math. J.} 18 (2008) 45--60; MR2467928.

\bibitem {Ekk1-tcs9}W. H. Ekkelkamp, The role of semismooth numbers in
factoring large numbers, \textit{Proc. Conf. on Algorithmic Number Theory},
ed. A.-M. Ernvall-Hyt\"{o}nen, M. Jutila, J. Karhum\"{a}ki and A. Lepist\"{o},
Turku Centre for Computer Science, 2007, pp. 40--44; http://oldtucs.abo.fi/publications/.

\bibitem {Ekk2-tcs9}W. H. Ekkelkamp, \textit{On the Amount of Sieving in
Factorization Methods}, Ph.D. thesis, Univ. Leiden, 2010; http://www.universiteitleiden.nl/en/research/research-output/.

\bibitem {BS-tcs9}E. Bach and J. Sorenson, Approximately counting semismooth
integers, \textit{Proc. 38th Internat. Symp. on Symbolic and Algebraic
Computation (ISSAC)}, ACM, 2013, pp. 23--30; arXiv:1301.5293; MR3206336.

\bibitem {Clff-tcs9}E. H. Cliffe, \textit{Reflections on the Number Field
Sieve}, Ph.D. thesis, Univ. of Bath, 2007; http://researchportal.bath.ac.uk/en/studentTheses/.

\bibitem {Trmr-tcs9}E. Tromer, \textit{Hardware-Based Cryptanalysis}, Ph.D.
thesis, Weizmann Institute of Science, 2007;
http://www.cs.tau.ac.il/\symbol{126}tromer/phd-dissertation/.

\bibitem {Grn1-tcs9}A. Granville, The anatomy of integers and permutations,
unpublished note, 2008, http://dms.umontreal.ca/\symbol{126}andrew/PDF/Anatomy.pdf.

\bibitem {Grn2-tcs9}A. Granville, J. Granville and R. J. Lewis, \textit{Prime
Suspects. The Anatomy of Integers and Permutations}, Princeton Univ. Press,
2019, pp. 200--201; MR3966460.

\bibitem {Llyd-tcs9}S. P. Lloyd, Ordered prime divisors of a random integer,
\textit{Annals of Probab.} 12 (1984) 1205--1212; MR0757777.

\bibitem {Kng2-tcs9}J. F. C. Kingman, The Poisson-Dirichlet distribution and
the frequency of large prime divisors, unpublished note, 2004, http://www.newton.ac.uk/documents/preprints/.

\bibitem {Tao-tcs9}T. Tao, Cycles of a random permutation, and irreducible
factors of a random polynomial, unpublished note, 2015, http://terrytao.wordpress.com/2015/07/15/.

\bibitem {AT-tcs9}R. Arratia and S. Tavar\'{e}, The cycle structure of random
permutations, \textit{Annals of Probab. }20 (1992) 1567--1591; MR1175278.

\bibitem {Fi4-tcs9}S. R. Finch, Golomb-Dickman constant, \textit{Mathematical
Constants}, Cambridge Univ. Press, 2003, pp. 284--292; MR2003519.

\bibitem {Fi5-tcs9}S. R. Finch, Extreme prime factors, \textit{Mathematical
Constants II}, Cambridge Univ. Press, 2019, pp. 171--172; MR3887550.

\bibitem {Fi6-tcs9}S. R. Finch, Components and cycles of random mappings,
\textit{forthcoming}.%

\begin{tabular}
[c]{lll}
& Steven Finch & \\
& MIT Sloan School of Management & \\
& Cambridge, MA, USA & \\
& \textit{steven\_finch@harvard.edu} &
\end{tabular}

\end{thebibliography}
\end{document}